\documentstyle[12pt]{article}
   \title{{\bf  Vertex operator algebras and the zeta function}}
    \author{J. Lepowsky}
    \date{\it Dedicated to Howard Garland on the occasion of his
sixtieth birthday}
    \begin{document}

    \bibliographystyle{alpha} 
    \maketitle

\begin{abstract}
We announce a new type of ``Jacobi identity'' for vertex operator
algebras, incorporating values of the Riemann zeta function at
negative integers.  Using this we ``explain'' and generalize some
recent work of S.~Bloch's relating values of the zeta function with
the commutators of certain operators and Lie algebras of differential
operators.
\end{abstract}

    \input amssym.def
    \input amssym  
 
    \newtheorem{rema}{Remark}[section]
    \newtheorem{propo}[rema]{Proposition}
    \newtheorem{theo}[rema]{Theorem}
    \newtheorem{defi}[rema]{Definition}
    \newtheorem{lemma}[rema]{Lemma}
    \newtheorem{corol}[rema]{Corollary}
    \newtheorem{exam}[rema]{Example}

\newcommand{\nordplus}{\mbox{\scriptsize ${+ \atop +}$}}
\newcommand{\nordbullet}{\mbox{\tiny ${\bullet\atop\bullet}$}}

\setcounter{section}{0}
\renewcommand{\theequation}{\thesection.\arabic{equation}}
\renewcommand{\therema}{\thesection.\arabic{rema}}
\setcounter{equation}{0}
\setcounter{rema}{0}

\section{Introduction}

Consider the famous classical ``formula''
\begin{equation}\label{1+}
1 + 2 + 3 + \cdots = - \frac{1}{12},
\end{equation}
which has the rigorous meaning
\begin{equation}\label{zeta-1}
\zeta(-1) = - \frac{1}{12}.
\end{equation}
Here $\zeta$ is of course the Riemann zeta function
\begin{equation}\label{zeta}
\zeta(s) = \sum_{n > 0} n^{-s}
\end{equation}
(analytically continued), and (\ref{1+}) is classically generalized by
the formal equality
\begin{equation}\label{zeta-s}
\sum_{n>0}n^s = \zeta(-s)
\end{equation}
for $s = 0,1,2,\dots$.  The classical number theory behind this
analytic continuation is well known to be related to the
widely-pervasive issue of regularizing certain infinities in quantum
field theory, in particular, in conformal field theory.  Here we shall
announce some general principles of vertex operator algebra theory
that elucidate the passage {}from the unrigorous but suggestive formula
(\ref{1+}) to formula (\ref{zeta-1}), and the generalization
(\ref{zeta-s}).  In the process, we shall ``explain'' some recent work
of S. Bloch's involving zeta-values and differential operators.  The
work \cite{L2} contains details and related results.  The material
that we shall present involves foundational notions of vertex operator
algebra theory, and we shall try to make this writeup accessible to
nonspecialists by reviewing elementary matters.

We were motivated by a desire to understand some very interesting
phenomena found by Bloch \cite{Bl} relating the values $\zeta (-n)$,
$n = 1, 3, 5, \dots$, of the zeta function at negative odd integers to
the commuatators of certain operators on an infinite-dimensional
space.  We shall begin with some elementary background and a brief
description of this work, then we shall explain how to recover and
somewhat generalize these results using vertex operator algebra
theory, and finally, we shall place these ideas and results into a
very general context and present some new general results in vertex
operator algebra theory.  These methods serve incidentally to enhance
the many already-existing motivations for vertex operator algebra
theory (see \cite{Bo}, \cite{FLM}) and its underlying formal calculus
(as developed in \cite{FLM} and \cite{FHL}).

One of our main themes is to ``always'' use generating functions---to
introduce new formal variables and generating functions in order to
try to make complicated things easier and more natural and at the same
time, much more general, as in the corresponding parts of \cite{FLM}.
We use commuting formal variables rather than complex variables
because they provide the most natural way to handle the
doubly-infinite series such as $\delta (x) = \sum_{n \in {\Bbb Z}}x^n$
that pervade the natural formulations and proofs.  Other central
themes are to exploit the formal exponential of the differential
operator $x{\frac{d}{dx}}$ as a formal change-of-variables
automorphism (again as in \cite{FLM}); to formulate Euler's
interpretation of the divergent series (\ref{zeta-s}) by means of the
operator product expansion in conformal field theory; and to place
considerations about Lie algebras of differential operators into the
very general context of what we termed the ``Jacobi identity''
\cite{FLM} for vertex operator algebras.  There are some interesting
points of contact between the present work and \cite{KR}, \cite{M} and
\cite{FKRW}.

I am very pleased to dedicate this paper to Howard Garland on the
occasion of his sixtieth birthday.  This work was presented in a talk
at Yale University in fall, 1997 at a seminar in his honor.  I would
like to mention here that it was {}from Howard Garland that Robert
Wilson and I, in 1980, first learned about the idea of using formal
delta-function calculus, which was also used in
\cite{DKM}; cf. \cite{FLM}.  This was one among many of Howard's 
insights that have influenced us.

This work was also presented in a talk at the May, 1998 Conference on
Representations of Affine and Quantum Affine Algebras and Related
Topics at North Carolina State University.  I would like to thank
Naihuan Jing and Kailash Misra, the organizers, for a stimulating
conference.  

I am very grateful to Spencer Bloch for informing me about his work
and for many valuable discussions.

This work was partially supported by NSF grants DMS-9401851 and
DMS-9701150.

\renewcommand{\theequation}{\thesection.\arabic{equation}}
\renewcommand{\therema}{\thesection.\arabic{rema}}
\setcounter{equation}{0}
\setcounter{rema}{0}

\section{Background}

Consider the commutative associative algebra $\Bbb C [t,t^{-1}]$ of
Laurent polynomials in an indeterminate $t$, and consider its Lie
algebra $\frak {d}$ of derivations:
\begin{equation}
{\frak {d}} = \mbox{\rm Der}\ \Bbb C [t,t^{-1}],
\end{equation}
the Lie algebra of formal vector fields on the circle, with basis
$\{t^n D | n \in \Bbb Z\}$, where
\begin{equation}\label{D}
D = D_t = t \frac{d}{dt}.
\end{equation}
(A preview of one of our main themes: The ``homogeneous'' differential
operator $D_t$, rather than $\frac{d}{dt}$, will be the appropriate
form of differentiation for our considerations, and we shall be using
it for various variables as well as $t$.)

Consider also the Virasoro algebra ${\frak {v}}$, the well-known central
extension
\begin{equation}\label{Vir}
0 \rightarrow \Bbb C c \rightarrow {\frak {v}} \rightarrow {\frak {d}}
\rightarrow 0,
\end{equation}
where ${\frak {v}}$ has basis $\{L(n) | n \in \Bbb Z \}$ together with
a central element $c$; the bracket relations among the $L(n)$ are
given by
\begin{equation}\label{Virbrackets}
[L(m),L(n)] = (m-n)L(m+n) + \frac{1}{12} (m^3 - m) \delta_{m+n,0} c
\end{equation}
and $L(n)$ maps to $-t^n D$ in (\ref{Vir}).  The number $\frac{1}{12}
(m^3 - m)$, being one-half the binomial coefficient ${m+1} \choose 3$,
is closely related to a third derivative, which becomes visible when
we use generating functions to write the bracket relations
(\ref{Virbrackets}), as we review below.  The Virasoro algebra is
naturally $\Bbb Z$-graded, with ${\rm deg}\ L(n) = n$ and ${\rm deg}\
c = 0$.

The following classical realization of the Lie algebra ${\frak {v}}$
is well known: We start with the Heisenberg Lie algebra with basis
consisting of the symbols $h(n)$ for $n \in \Bbb Z$, $n \neq 0$ and a
central element $1$, with the bracket relations
\begin{equation}\label{Heisbrackets}
[h(m),h(n)] = m \delta_{m+n,0} 1.
\end{equation}
For convenience we adjoin an additional central basis element $h(0)$,
so that the relations (\ref{Heisbrackets}) hold for all $m,n \in
\Bbb Z$.  This Lie algebra acts irreducibly on the polynomial algebra
\begin{equation}\label{S}
S = \Bbb C [h(-1),h(-2),h(-3),\dots]
\end{equation}
as follows: For $n<0$, $h(n)$ acts as the multiplication operator; for
$n>0$, $h(n)$ acts as the operator $n \frac{\partial}{\partial
h(-n)}$; $h(0)$ acts as 0; and $1$ acts as the identity operator.
Then ${\frak {v}}$ acts on $S$ by means of the following operators:
\begin{equation}\label{c}
c \mapsto 1,
\end{equation}
\begin{equation}\label{Ln}
L(n) \mapsto \frac{1}{2} \sum_{j \in \Bbb Z} h(j)h(n-j) \ \ \mbox{for}
\ n \neq 0,
\end{equation}
\begin{equation}\label{L0}
L(0) \mapsto \frac{1}{2} \sum_{j \in \Bbb Z} h(-|j|)h(|j|).
\end{equation}

It is important to observe that in the case of $L(0)$, the absolute
values make the formal quadratic operator well defined, while for $n
\neq 0$ the operator is well defined as it stands, since
$[h(j),h(n-j)] = 0$.  Thus the operators (\ref{Ln}) and (\ref{L0}) are
in ``normal-ordered form,'' that is, the ``annihilation operators''
$h(n)$ for $n>0$ act to the right of the ``creation operators'' $h(n)$
for $n<0$.  Using colons to denote normal ordering (the appropriate
reordering of the factors in the product if necessary), we thus have
\begin{equation}\label{Lnnormal}
L(n) \mapsto \frac{1}{2} \nordbullet\sum_{j \in \Bbb Z}
h(j)h(n-j)\nordbullet
\end{equation}
for all $n \in \Bbb Z$.

It is an instructive and not-so-trivial (classical) exercise to verify
by direct computation that the operators (\ref{Lnnormal}) indeed
satisfy the bracket relations (\ref{Virbrackets}).  (This exercise and
the related constructions are presented in \cite{FLM}, for example,
where the standard generalization of this construction of ${\frak
{v}}$ using a Heisenberg algebra based on a finite-dimensional space
of operators $h(n)$ for each $n$ is also carried out.)

Vertex operator algebra theory and conformal field theory place this
exercise into a very general, natural setting (among many other
things), with conceptual approaches and techniques (cf. \cite{FLM}).
It is standard procedure to embed operators such as $h(n)$ and $L(n)$
into generating functions and to compute with these generating
functions, using a formal calculus, and to systematically avoid
computing with the individual operators.  Doing this vastly simplifies
computations that would otherwise be complicated or sometimes almost
impossible.  In fact, we shall be using a number of
generating-function ideas below.

The space $S$ carries a natural $\Bbb Z$-grading, determined by the
rule ${\rm deg}\ h(j) = j$ for $j < 0$.  Then $S$ is in fact graded by
the nonpositive integers, and the ${\frak {v}}$-module $S$ is a graded
module.  It turns out to be appropriate to use the negative of this
grading, that is, to define a new grading (by ``conformal weights'')
on the space $S$ by the rule ${\rm wt}\ h(-j) = j$ for $j > 0$.  One
reason why this is natural is that for each $n \geq 0$, the
homogeneous subspace of $S$ of weight $n$ coincides with the
eigenspace of the operator $L(0)$ with eigenvalue $n$, as is easy to
see.  For $n \in {\Bbb Z}$ (or $n \geq 0$) we define $S_n$ to be the
homogeneous subspace of $S$ of weight $n$, and we consider the formal
power series in the formal variable $q$ given by
\begin{equation}
\mbox{dim}_* S = \sum_{n \geq 0} (\mbox{dim}\;S_n) q^n
\end{equation}
(the ``graded dimension'' of the graded space $S$).  Clearly, {}from the
definitions,
\begin{equation}
\mbox{dim}_* S = \prod_{n > 0} (1 - q^n)^{-1}.
\end{equation}

Here are the main points about these classical considerations that we
want to emphasize: As is well known in conformal field theory,
removing the normal ordering in the definition of the operator $L(0)$
introduces an infinity which formally equals $\frac{1}{2} \zeta(-1)$:
The unrigorous expression
\begin{equation}\label{Lbar(0)}
{\bar L}(0) = \frac{1}{2} \sum_{j \in \Bbb Z} h(-j)h(j)
\end{equation}
formally equals (by the bracket relations (\ref{Heisbrackets}))
\begin{equation}\label{Lbar(0)2}
L(0) + \frac{1}{2} (1 + 2 + 3 + \cdots),
\end{equation}
which itself formally equals
\begin{equation}\label{Lbar(0)3}
L(0) + \frac{1}{2} \zeta(-1) = L(0) - \frac{1}{24}.
\end{equation}
We rigorize ${\bar L}(0)$ by defining it as:
\begin{equation}\label{Lbar(0)rig}
{\bar L}(0) = L(0) + \frac{1}{2} \zeta(-1),
\end{equation}
and we define
\begin{equation}\label{Lbar(n)}
{\bar L}(n) = L(n) \ \ \mbox{for} \ n \neq 0,
\end{equation}
to get a new basis of the Lie algebra ${\frak {v}}$.  (We are
identifying the elements of ${\frak {v}}$ with operators on the space
$S$.)  The brackets become:
\begin{equation}\label{newVirbrackets}
[{\bar L}(m),{\bar L}(n)] = (m-n){\bar L}(m+n) + \frac{1}{12} m^3
\delta_{m+n,0};
\end{equation}
that is, $m^3 - m$ in (\ref{Virbrackets}) has become the pure monomial
$m^3$.

Perhaps the most important (and also well-known) thing accomplished by
this formal removal of the normal ordering is the following: We define
a new grading of the space $S$ by using the eigenvalues of the
modified operator ${\bar L}(0)$ in place of $L(0)$, so that the
grading of $S$ is ``shifted'' {}from the previous grading by conformal
weights by the subtraction of $\frac{1}{24}$ {}from the weights.  We let
$\chi (S)$ be the corresponding graded dimension, so that
\begin{equation}\label{chiS}
\chi (S) = \frac{1}{\eta (q)},
\end{equation}
where
\begin{equation}
\eta (q) = q^{\frac{1}{24}} \prod_{n > 0} (1 - q^n). 
\end{equation}
The point is that $\eta (q)$ has important (classical) modular
transformation properties, unlike $\prod_{n > 0} (1 - q^n)$, when
viewed as a function of $\tau$ in the upper half-plane via the
substitution $q = e^{2 \pi i \tau}$; $\eta (q)$ is Dedekind's
eta-function.

Bloch \cite{Bl} extended this classical story in various ways, in
particular, the following: Instead of the Lie algebra $\frak {d}$, we
consider the larger Lie algbebra of formal differential operators,
spanned by
\begin{equation}
\{ t^n D^m | n \in \Bbb Z, \ m \geq 0 \}
\end{equation}
or more precisely, we restrict to $m > 0$ and further, to the Lie
subalgebra ${\cal D}^+$ spanned by the differential operators of the
form $D^r (t^nD) D^r$ for $r \geq 0, \ n \in \Bbb Z$.  Then we can
construct a central extension of ${\cal D}^+$ using generalizations of
the normally-ordered quadratic operators (\ref{Lnnormal}) above:
\begin{equation}\label{Lrn}
L^{(r)}(n) = \frac{1}{2} \sum_{j \in \Bbb Z} j^r h(j) (n-j)^r h(n-j) \
\ \mbox{for} \ n \neq 0,
\end{equation}
\begin{equation}\label{Lr0}
L^{(r)}(0) = \frac{1}{2} \sum_{j \in \Bbb Z} (-j)^r h(-|j|) j^r
h(|j|),
\end{equation}
that is,
\begin{equation}\label{Lrnnormal}
L^{(r)}(n) = \nordbullet\frac{1}{2} \sum_{j \in \Bbb Z} j^r h(j) (n-j)^r
h(n-j)\nordbullet
\end{equation}
for $n \in \Bbb Z$.  Direct computation of the commutators among these
operators \cite{Bl} shows that they provide a central extension of
${\cal D}^+$ such that
\begin{equation}
L^{(r)}(n) \mapsto (-1)^{r+1} D^r (t^nD) D^r
\end{equation}
(cf. \cite{KP}).  (It is not surprising in retrospect that these
operators $L^{(r)}(n)$ are related to differential operators, because
the generating function of these operators as $n$ ranges through $\Bbb
Z$ is based on $D^r$, as we discuss below.)

A central point of \cite{Bl} is that the formal removal of the
normal-ordering procedure in the definition (\ref{Lr0}) of
$L^{(r)}(0)$ adds the infinity $(-1)^r \frac{1}{2} \zeta(-2r-1) =$ ``
$\sum_{n>0} n^{2r+1}$'' (generalizing
(\ref{Lbar(0)})--(\ref{Lbar(0)3})), and if we correspondingly define
\begin{equation}\label{lbarr(0)}
{\bar L}^{(r)} (0) = L^{(r)} (0) + (-1)^r \frac{1}{2} \zeta(-2r-1)
\end{equation}
and ${\bar L}^{(r)} (n) = L^{(r)} (n)$ for $n \neq 0$ (generalizing
(\ref{Lbar(0)rig}) and (\ref{Lbar(n)})), the commutators simplify in a
remarkable way: As direct computation \cite{Bl} shows, the complicated
polynomial in the scalar term of $[{\bar L}^{(r)} (m),{\bar L}^{(s)}
(-m)]$ reduces to a pure monomial in $m$, by analogy with, and
generalizing, the passage {}from $m^3 - m$ to $m^3$ in
(\ref{newVirbrackets}).  The precise formulas can be found in
\cite{Bl}, along with further results; for instance, in \cite{Bl},
these considerations and results are generalized to Dirichlet
$L$-series in place of the zeta function.

\renewcommand{\theequation}{\thesection.\arabic{equation}}
\renewcommand{\therema}{\thesection.\arabic{rema}}
\setcounter{equation}{0}
\setcounter{rema}{0}

\section{First ``explanation'' and generalization}

Our goal is to present two layers of ``explanation'' and
generalization of the results of \cite{Bl} sketched above.  First we
need some elementary formal background:

What does $\zeta(-2r-1)$ ``mean,'' for a nonnegative integer $r$?

It is a well-known classical fact that for $k > 1$,
\begin{equation}\label{zetaB}
\zeta(-k+1) = - \frac{B_{k}}{k},
\end{equation}
where the $B_{k}$ are the Bernoulli numbers, defined by the generating
function
\begin{equation}\label{Ber}
\frac{x}{e^x - 1} = \sum_{k \geq 0} \frac{B_{k}}{k!} x^k,
\end{equation}
where $x$ is a formal variable.  This formal power series in $x$ is
understood to be computed (on the left-hand side) by expanding $e^x -
1$ as the formal series $x + \frac{x^2}{2!} + \cdots$ and performing
the division of formal power series to obtain a formal power series
with constant term $1$; this of course corresponds to expanding a
complex function in a certain domain, but we are operating purely
formally.

Why does $B_k$ defined in this way have anything to do with the formal
series $\sum_{n>0} n^{k-1}$?  We recall Euler's heuristic
interpretation of such formal sums as ``$1 + 2 + \cdots$'' (cf. the
Preface of \cite{Hi}); actually, we give a variant of Euler's
interpretation adapted to the main theme that we shall introduce:

Consider the expansion ({}from (\ref{Ber}))
\begin{equation}\label{fromBer}
\frac{1}{1 - e^x} = - \sum_{k \geq 0} \frac{B_{k}}{k!} x^{k-1}.
\end{equation}
Expand the left-hand side {\it unrigorously} as the formal geometric
series
\begin{equation}\label{unrig}
1 + e^x + e^{2x} + \cdots = 1 + \sum_{k \geq 0} \frac{1^k}{k!} x^k +
\sum_{k \geq 0} \frac{2^k}{k!} x^k + \cdots.
\end{equation}
For $k>1$, the coefficient of $x^{k-1}$ in this formal expression is
\begin{equation}
\frac{1}{(k-1)!} (1^{k-1} + 2^{k-1} + \cdots),
\end{equation}
which looks like $\frac{1}{(k-1)!} \zeta (-k+1)$.  Also, the
coefficient of $x^0$ in (\ref{unrig}) is formally $1 + \frac{1}{0!}
(1^0 + 2^0 + \cdots)$, which we formally view as $1 + \zeta (0)$ (and
not as $\zeta (0)$).  Thus, formally equating the coefficients of
$x^l$ for $l \geq 0$ in (\ref{Ber}) ``explains'' (\ref{zetaB}) and the
fact that $\zeta(0) = - B_1 - 1$ ($= - \frac{1}{2}$); now we know what
(\ref{zeta-s}) says.

The key point here is the interplay between the formal geometric
series expansion (in powers of $e^x$) and the expansion in powers of
$x$.

Now, how do we interpret all of this via vertex operator algebra
theory?

First note that the expressions (\ref{Lrn})--(\ref{Lrnnormal}) above
for $L^{(r)}(n)$ suggest $r^{\rm th}$ derivatives.  We have already
mentioned that a basic theme in vertex operator algebra theory is to
always use appropriate generating functions (as we just did,
incidentally, in the heuristic discussion above).  First we put our
individual operators into generating functions.  Using a formal
variable $x$, we define
\begin{equation}\label{h(x)}
h(x) = \sum_{n \in \Bbb Z} h(n)x^{-n}
\end{equation}
and
\begin{equation}
L^{(r)}(x) = \sum_{n \in \Bbb Z} L^{(r)}(n)x^{-n},
\end{equation}
and using $D_x$ to denote the operator $x \frac{d}{dx}$ (recall the
comment after (\ref{D})), we observe that
\begin{equation}\label{Lrx}
L^{(r)}(x) = {\frac{1}{2}}\nordbullet(D_x^r h(x))^2\nordbullet,
\end{equation}
where the colons, as always, denote normal ordering (recall
(\ref{Lnnormal})).  (For other purposes, other versions of these
generating functions are used, in particular, $h(x) = \sum_{n \in \Bbb
Z} h(n)x^{-n-1}$, as in (\ref{Yh(-1)}) below, in place of
(\ref{h(x)}), but we have chosen the appropriate generating functions
for our purposes.)

Using standard elementary techniques, we could directly compute the
brackets $[L^{(r)}(x_1),L^{(s)}(x_2)]$ of these generating functions,
for $r,s \ge 0$, where $x_1$ and $x_2$ are independent commuting
formal variables.  (As always in vertex operator algebra theory or
conformal field theory, when we consider such operations as brackets
of generating functions, we need independent commuting formal
variables; the expression $[L^{(r)}(x),L^{(s)}(x)]$, with the variable
$x$ repeated, would be meaningless.)  But this computation, which
might be carried out as a more complicated variant of the argument on
pp. 224--226 of \cite{FLM}, for example, would not be simple.  It
would of course recover the information of the brackets
$[L^{(r)}(m),L^{(s)}(n)]$ computed in \cite{Bl}.

The best use of generating functions in this context is instead to
also introduce suitable generating functions over the number of {\it
derivatives}.  Consider the elementary formal Taylor theorem equating
the application of a formal exponential of a formal multiple of
$\frac{d}{dx}$ with a formal substitution operation:
\begin{equation}\label{Taylor}
e^{y \frac{d}{dx}} f(x) = f(x+y),
\end{equation}
where $f(x)$ is an arbitrary formal series of the form $\sum_n a_n
x^n$, and where it is understood that each binomial expression
$(x+y)^n$ is to be expanded in nonnegative integral powers of $y$.
Here $n$ is allowed to range over something very general, like $\Bbb
Z$ or even $\Bbb C$, say, and the $a_n$ lie in a fixed vector space;
$f(x)$ is very definitely not necessarily the expansion of an analytic
function.  Formula (\ref{Taylor}) is proved by direct formal expansion
of both sides (cf. \cite{FLM}, Proposition 8.3.1; Taylor's theorem in
this kind of generality is heavily exploited in Chapter 8 of
\cite{FLM}, for instance).  Now $\frac{d}{dx}$ is of course a formal
infinitesimal translation (as (\ref{Taylor}) states), but for our
present purposes we want the following formal multiplicative analogue
of (\ref{Taylor}):
\begin{equation}\label{infinitdil}
e^{y D_x} f(x) = f(e^y x),
\end{equation}
with $f(x)$ as above (again cf. \cite{FLM}, Proposition 8.3.1), which
expresses the fact that $D_x$ is a formal infinitesimal dilation.

Now $\nordbullet(D_x^r h(x))^2\nordbullet$ (recall (\ref{Lrx})) is
hard to put into a ``good'' generating function over $r$, but we can
make the problem easier by making it more general: Consider
independently many derivatives on each of the two factors $h(x)$ in
$\nordbullet h(x)^2\nordbullet$, use two new independent formal
variables $y_1$ and $y_2$, and form the generating function
\begin{equation}\label{Ly1y2}
L^{(y_1,y_2)}(x) = {\frac{1}{2}}
\nordbullet(e^{y_1 D_x} h(x))(e^{y_2 D_x}
h(x))\nordbullet = {\frac{1}{2}}
\nordbullet h(e^{y_1} x)h(e^{y_2} x)\nordbullet
\end{equation}
(where we use (\ref{infinitdil})), so that $L^{(r)}(x)$ is a
``diagonal piece'' of this generating function in the sense that it is
$(r!)^2$ times the coefficient of $y_1^r y_2^r$ in $L^{(y_1,y_2)}(x)$.
Using formal vertex operator calculus techniques (generalizing the
argument on pp. 224-226 of \cite{FLM}, for example), we can calculate
\begin{equation}\label{bracketofquadratics}
[\nordbullet h(e^{y_1} x_1)h(e^{y_2} x_1)\nordbullet,\nordbullet
h(e^{y_3} x_2)h(e^{y_4} x_2)\nordbullet].
\end{equation}

Then, a nontrivial, and in fact quite tricky, vertex operator
computation recovers Bloch's formulas, in somewhat generalized form,
as we explain next.  Here are the main points:

The expression $h(x)^2$ is not rigorous (as we observe for instance by
trying to compute the coefficient of any fixed power of $x$ in this
expression), while the normal-ordered expression $\nordbullet
h(x)^2\nordbullet$ is of course rigorous.  The expression $h(e^{y_1}
x)h(e^{y_2} x)$ is still not rigorous (even though the expressions
$e^{y_1} x$ and $e^{y_2} x$ are distinct), as we see by (for example)
trying to compute the constant term in the variables $y_1$ and $y_2$
in this expression.  The lack of rigor in this expression in fact
corresponds exactly to the occurrence of formal sums like $\sum_{n>0}
n^r$ with $r>0$, as we have been discussing above.

However, we have
\begin{equation}\label{hx1hx2}
h(x_1)h(x_2) = \nordbullet h(x_1)h(x_2)\nordbullet + x_2
\frac{\partial}{\partial x_2}
\frac{1}{1 - x_2 / x_1}
\end{equation}
(an exercise using elementary vertex operator techniques), and it
follows that
\begin{equation}\label{hex1hex2}
h(e^{y_1} x_1)h(e^{y_2} x_2) = \nordbullet h(e^{y_1} x_1)h(e^{y_2}
x_2)\nordbullet + x_2
\frac{\partial}{\partial x_2} \frac{1}{1 - e^{y_2} x_2 / e^{y_1}
x_1}.
\end{equation}
Note that $x_2\frac{\partial}{\partial x_2}$ can be replaced by $-
\frac{\partial}{\partial y_1}$ in the last expression (and this is one
illustration of the naturalness of our emphasis on the ``homogeneous''
differential operator $D_x = x\frac{\partial}{\partial x}$ rather than
$\frac{\partial}{\partial x}$).  The expression $\frac{1}{1 - e^{y_2}
x_2 / e^{y_1} x_1}$ came {}from, and is, a geometric series expansion
(recall (\ref{hx1hx2})).

Now we try to set $x_1 = x_2 \ (= x)$ in (\ref{hex1hex2}).  The result
of this procedure is unrigorous on the left-hand side, as we have
pointed out above, {\it but the result has rigorous meaning on the
right-hand side}, because the normal-ordered product $\nordbullet
h(e^{y_1} x)h(e^{y_2} x)\nordbullet$ is certainly well defined, and
the expression $-
\frac{\partial}{\partial y_1} \frac{1}{1 - e^{-y_1 + y_2}}$ can be
interpreted rigorously as in (\ref{Ber}) and (\ref{fromBer}); more
precisely (the role of $x$ in (\ref{fromBer}) being played here by
$-y_1 + y_2$), we take $\frac{1}{1 - e^{-y_1 + y_2}}$ to mean the
formal (Laurent) series in $y_1$ and $y_2$ of the shape
\begin{equation}\label{defof1/1-e}
\frac{1}{1 - e^{-y_1 + y_2}} = (y_1 - y_2)^{-1}F(y_1,y_2),
\end{equation}
where $(y_1 - y_2)^{-1}$ is understood as the binomial expansion
(geometric series) in nonnegative powers of $y_2$ and $F(y_1,y_2)$ is
an (obvious) formal power series in (nonnegative powers of) $y_1$ and
$y_2$.  This motivates us to define a new ``normal-ordering''
procedure
\begin{equation}\label{hexhex}
\nordplus h(e^{y_1} x)h(e^{y_2} x)\nordplus = 
\nordbullet h(e^{y_1} x)h(e^{y_2} x)\nordbullet -
\frac{\partial}{\partial y_1} \frac{1}{1 - e^{-y_1 + y_2}},
\end{equation}
with the last part of the right-hand side being understood as we just
indicated.  Again compare this with the heuristic discussion above;
this expression came {}from a geometric series, but it becomes rigorous
only when we expand in the new way (actually, we might alternatively
replace the binomial expansion $(y_1 - y_2)^{-1}$ by the different
expansion of the same formal expression in nonnegative powers of $y_1$
rather than of $y_2$, but it is more natural to make the choice that
we did).

Formula (\ref{hexhex}) and its indicated interpretation give a natural
``explanation'' of the zeta-function-modified operators defined in
(\ref{lbarr(0)}):  We use (\ref{hexhex}) to define the following
analogues of the operators (\ref{Ly1y2}):
\begin{equation}\label{Lbary1y2}
{\bar L}^{(y_1,y_2)}(x) = 
{\frac{1}{2}}
\nordplus h(e^{y_1} x)h(e^{y_2} x)\nordplus,
\end{equation}
and it is easy to check that the operator ${\bar L}^{(r)} (n)$ is exactly
$(r!)^2$ times the coefficient of $y_1^r y_2^r x_0^{-n}$ in
(\ref{Lbary1y2}); the significant case is the case $n = 0$.

We are now ready to formulate the result mentioned just after
(\ref{bracketofquadratics}) above.  With the new normal ordering
(\ref{hexhex}) replacing the old one, remarkable cancellation occurs
in the commutator (\ref{bracketofquadratics}), and the result is:
\begin{theo}\label{theoremforLbar}
With the formal delta-function Laurent series $\delta(x)$ defined as
\begin{equation}\label{delta}
\delta(x) = \sum_{n \in {\Bbb Z}} x^n,
\end{equation}
and with independent commuting formal variables as indicated, we have:
\begin{eqnarray}\label{Lbarbrackets}
\lefteqn{[{\bar L}^{(y_1,y_2)}(x_1),{\bar L}^{(y_3,y_4)}(x_2)]}\nonumber\\
&&= - {\frac{1}{2}} \frac{\partial}{\partial y_1}	
\biggl({\bar L}^{(-y_1+y_2+y_3,y_4)}(x_2)
\delta \left({\frac{e^{y_1}x_1}{e^{y_3}x_2}}\right)\nonumber\\
&&\hspace{2em} + {\bar L}^{(-y_1+y_2+y_4,y_3)}(x_2)
\delta \left({\frac{e^{y_1}x_1}{e^{y_4}x_2}}\right)\biggr)\nonumber\\
&&\hspace{2em} - {\frac{1}{2}} \frac{\partial}{\partial y_2}	
\biggl({\bar L}^{(y_1-y_2+y_3,y_4)}(x_2)
\delta \left({\frac{e^{y_2}x_1}{e^{y_3}x_2}}\right)\nonumber\\
&&\hspace{2em} + {\bar L}^{(y_1-y_2+y_4,y_3)}(x_2)
\delta \left({\frac{e^{y_2}x_1}{e^{y_4}x_2}}\right)\biggr).
\end{eqnarray}
\end{theo}

Actually, hidden in the right-hand side of (\ref{Lbarbrackets}) are
formal expressions (suitably expanded) of the type $(y_1 - y_2 - y_3 +
y_4)^{-3}$ because of the formal pole $y_1 = y_2$ in (\ref{hexhex}),
and these expressions, multiplied by the formal delta-function
expressions, are the source of the pure monomials in $m$ that we set
out to explain (recall the end of Section 2).  Indeed, the
delta-function expression $\delta
\left(e^{y_1}x_1 / e^{y_3}x_2 \right)$, for instance, can be
written as $e^{y_1 D_{x_1}}e^{y_3 D_{x_2}}\delta
\left(x_1 / x_2 \right)$, and when we extract and equate the
coefficients of the monomials in the variables $y_1^r y_2^r y_3^s
y_4^s$ on the two sides of (\ref{Lbarbrackets}), we get expressions
like $(D^j \delta)\left(x_1 / x_2 \right)$, whose expansion, in
turn, in powers of $x_1$ and $x_2$ clearly yields a pure monomial
analogous to and generalizing the expression $m^3$ in
(\ref{newVirbrackets}).  In this way, we have an explicit
generalization and explanation of Bloch's formula for $[{\bar L}^{(r)}
(m),{\bar L}^{(s)} (n)]$ in terms of the operators ${\bar L}^{(j)}
(m+n)$ and a monomial in $m$.

We have been working all along with a Heisenberg algebra based on a
one-dimensional space---that is, a Heisenberg algebra with only one
dimension of operators, spanned by the element $h(n)$, for each $n$,
but all of these considerations hold equally well in the more general
situation where we start with a Heisenberg algebra based on a
finite-dimensional space.

All of this has been an interesting (and nontrivial) special case of
vertex operator algebra theory, but what is {\it really} happening?

\renewcommand{\theequation}{\thesection.\arabic{equation}}
\renewcommand{\therema}{\thesection.\arabic{rema}}
\setcounter{equation}{0}
\setcounter{rema}{0}

\section{Second ``explanation'' and generalization}

It turns out that Theorem \ref{theoremforLbar} is an extremely special
case of a something that can be formulated and proved for an arbitrary
vertex operator algebra (and indeed this gives another motivation for
the general theory).  We recall the definition of the notion of vertex
(operator) algebra {}from \cite{Bo},
\cite{FLM} and \cite{FHL}; the principles that we have found are based
heavily on the ``Jacobi identity'' as formulated in \cite{FLM} and
\cite{FHL}:

\begin{defi}\label{VOA}
{\rm A {\it vertex operator algebra} $(V, Y, {\bf 1}, \omega)$, or simply $V$
(over ${\Bbb C}$), is a ${\Bbb
Z}$-graded vector space (graded by {\it weights})
\begin{equation}
V=\coprod_{n\in {\Bbb Z}}V_{(n)}; \ \mbox{\rm for}\ v\in V_{(n)},\;n=\mbox{\rm wt}\ v;
\end{equation}
such that
\begin{equation}
\mbox{\rm dim }V_{(n)}<\infty\;\;\mbox{\rm for}\; n \in {\Bbb Z},
\end{equation}
\begin{equation}
V_{(n)}=0\;\;\mbox{\rm for} \;n\; \mbox{\rm sufficiently small},
\end{equation}
equipped with a linear map  $V\otimes V\to V[[x, x^{-1}]]$, or
equivalently,
\begin{eqnarray}
V&\to&(\mbox{\rm End}\; V)[[x, x^{-1}]]\nonumber \\
v&\mapsto& Y(v, x)={\displaystyle \sum_{n\in{\Bbb Z}}}v_{n}x^{-n-1}
\;\;(\mbox{\rm where}\; v_{n}\in
\mbox{\rm End} \;V),
\end{eqnarray}
$Y(v, x)$ denoting the {\it vertex operator associated with} $v$, and
equipped also with two distinguished homogeneous vectors ${\bf 1}\in
V_{(0)}$ (the {\it vacuum}) and $\omega \in V_{(2)}$. The following
conditions are assumed for $u, v \in V$: the {\it lower truncation
condition} holds:
\begin{equation}
u_{n}v=0\;\;\mbox{\rm for}\;n\; \mbox{\rm sufficiently large}
\end{equation}
(or equivalently, $Y(u,x)v$ involves only finitely many negative
powers of $x$);
\begin{equation}
Y({\bf 1}, x)=1\;\; (1\;\mbox{\rm on the right being the identity
operator});
\end{equation}
the {\it creation property} holds:
\begin{equation}
Y(v, x){\bf 1} \in V[[x]]\;\;\mbox{\rm and}\;\;\lim_{x\rightarrow
0}Y(v, x){\bf 1}=v
\end{equation}
(that is, $Y(v, x){\bf 1}$ involves only nonnegative integral powers
of $x$ and the constant term is $v$); with $\delta(x)$ as in
(\ref{delta}) and with binomial expressions understood (as above) to
be expanded in nonnegative powers of the second variable, the {\it
Jacobi identity} (the main axiom) holds:
\begin{eqnarray}\label{jacobi}
&x_{0}^{-1}\delta
\left({\displaystyle\frac{x_{1}-x_{2}}{x_{0}}}\right)Y(u, x_{1})Y(v,
x_{2})-x_{0}^{-1} \delta
\left({\displaystyle\frac{x_{2}-x_{1}}{-x_{0}}}\right)Y(v, x_{2})Y(u,
x_{1})&\nonumber \\ &=x_{2}^{-1} \delta
\left({\displaystyle\frac{x_{1}-x_{0}}{x_{2}}}\right)Y(Y(u, x_{0})v,
x_{2})&
\end{eqnarray}
(note that when each expression in (\ref{jacobi}) is applied to any
element of $V$, the coefficient of each monomial in the formal
variables is a finite sum; on the right-hand side, the notation
$Y(\cdot, x_{2})$ is understood to be extended in the obvious way to
$V[[x_{0}, x^{-1}_{0}]]$); the Virasoro algebra relations hold (acting
on $V$):
\begin{equation}
[L(m), L(n)]=(m-n)L(m+n)+{\displaystyle\frac{1}{12}}
(m^{3}-m)\delta_{n+m,0}({\rm rank}\;V)1
\end{equation}
for $m, n \in {\Bbb Z}$, where
\begin{equation}
L(n)=\omega _{n+1}\;\; \mbox{\rm for} \;n\in{\Bbb Z}, \;\;{\rm
i.e.},\;\;Y(\omega, x)=\sum_{n\in{\Bbb Z}}L(n)x^{-n-2}
\end{equation}
and 
\begin{eqnarray}
&{\rm rank}\;V\in {\Bbb C};&\\ 
&L(0)v=nv=(\mbox{\rm wt}\ v)v\;\;\mbox{\rm for}\;n \in {\Bbb
Z}\;\mbox{\rm and}\;v\in V_{(n)};&\\ 
&{\displaystyle \frac{d}{dx}}Y(v,x)=Y(L(-1)v, x)&
\end{eqnarray}
(the {\it  $L(-1)$-derivative property}).}
\end{defi}

Note that as in Theorem \ref{theoremforLbar}, the Jacobi identity
encapsulates higher derivatives of delta-function expressions, since
the expression $\delta ((x_1 - x_2)/x_0)$, say, can be expanded by
means of (\ref{Taylor}).  The use of the three formal variables,
rather than complex variables (which would require changes in the
formulas), allows the full symmetry of the Jacobi identity to reveal
itself (see \cite{FLM} and \cite{FHL}).

The polynomial algebra $S$ (recall (\ref{S})) carries a canonical
vertex operator algebra structure of rank 1 with vacuum vector ${\bf
1}$ equal to $1 \in S$ and with the operators $L(n)$ agreeing with the
operators defined in Section 2 (cf. \cite{FLM}).  We will not describe
the vertex operators $Y(v,x)$ here, except to say that
\begin{equation}\label{Yh(-1)}
Y(h(-1),x) = x^{-1}h(x) = \sum_{n \in \Bbb Z} h(n)x^{-n-1}
\end{equation}
(recall (\ref{h(x)})) and that the element $\omega$ is $\frac{1}{2}
(h(-1))^{2} \in S$.

If we omit the grading and the Virasoro algebra {}from Definition
\ref{VOA} and adjust the axioms appropriately, we have the notion of
``vertex algebra'' as introduced in \cite{Bo}, except that Borcherds
used certain special cases of the Jacobi identity instead of
(\ref{jacobi}).  The identity (\ref{jacobi}) is the canonical
``maximal'' axiom: It contains the ``full'' necessary information in
compact form; it is analogous to the classical Jacobi identity in the
definition of the notion of Lie algebra; and it is invariant in a
natural sense under the symmetric group on three letters (see
\cite{FLM} and \cite{FHL}).

There are also ``minimal'' axioms, stemming {}from the fact that the
(suitably formulated) ``commutativity'' of the operators $Y(u,x_1)$
and $Y(v,x_2)$ implies ``associativity'' (again suitably formulated)
and hence the Jacobi identity (see \cite{FLM} and \cite{FHL}; cf.
\cite{BPZ} and \cite{G}).  The simplest ``minimal'' axiom, as found 
in \cite{DL} (actually, in the greater generality of ``abelian
intertwining algebras'') states that for $u,v
\in V$, there exists $n \ge 0$ such that
\begin{equation}\label{weakcomm}
(x_1 - x_2)^n [Y(u,x_1),Y(v,x_2)] = 0
\end{equation}
(see \cite{DL}, formula (1.4)).  However, it is still a nontrivial
matter to construct examples, even relatively simple ones like $S$, of
vertex operator algebras, partly because one has to extend the
condition (\ref{weakcomm}) {}from generators to arbitrary elements of
$V$.  A general and systematic approach and solution to this and
related problems was first found by Li (see \cite{Li1}, \cite{Li2}).

The commutativity condition asserts that for $u,v \in V$,
\begin{equation}
Y(u,x_1)Y(v,x_2) \sim Y(v,x_2)Y(u,x_1),
\end{equation}
where ``$\sim$'' denotes equality up to a suitable kind of generalized
analytic continuation, and the associativity condition asserts that
\begin{equation}\label{assoc}
Y(u,x_1)Y(v,x_2) \sim Y(Y(u,x_1 - x_2)v,x_2),
\end{equation}
where the right-hand side and the generalized analytic continuation
have to be understood in suitable ways (again see \cite{FLM} and
\cite{FHL} and cf. \cite{BPZ} and \cite{G}); the right-hand side of 
(\ref{assoc}) is {\it not} a well-defined formal series in $x_1$ and
$x_2$.

On the level of these basic principles, for any vertex operator
algebra $V$ we shall now conceptually formulate and considerably
generalize the normal-ordering procedure (\ref{hexhex}) and we shall
formulate a new general ``Jacobi identity'' which implies Theorem
\ref{theoremforLbar} in the very particular case of the vertex
operator algebra $S$ and very special vertex operators.  Ideas in
Zhu's work \cite{Z1}, \cite{Z2} enter into our considerations.

Formally replacing $x_1$ by $e^{y}x_2$ in (\ref{assoc}), we find
(formally and unrigorously) that 
\begin{equation}\label{assoc2}
Y(u,e^{y}x_2)Y(v,x_2) \sim Y(Y(u,(e^{y}-1)x_2)v,x_2).
\end{equation}
Now we observe that while the left-hand side of (\ref{assoc2}) is not
a well-defined formal series in the formal variables $y$ and $x_2$,
the right-hand side of (\ref{assoc2}) {\it is} in fact a well-defined
formal series in these formal variables.  By replacing $x_1$ by
$e^{y}x_2$ we have made the right-hand side of (\ref{assoc2}) rigorous
(and the left-hand side unrigorous).  This situation should be
compared with our motivation for introducing the normal-ordering
procedure (\ref{hexhex}) above.

Next, instead of the vertex operators $Y(v,x)$, we want the modified
vertex operators defined for homogeneous elements $v
\in V$ by:
\begin{equation}
X(v,x) = x^{{\rm wt}\;v}Y(v,x) = Y(x^{L(0)}v,x),
\end{equation}
as in \cite{FLM}, formula (8.5.27) (recall that $L(0)$-eigenvalues
define the grading of $V$); the formula $X(v,x) = Y(x^{L(0)}v,x)$
works for {\it all} $v \in V$ (not necessarily homogeneous).  For
instance, the operator $h(x)$ (\ref{h(x)}) is exactly $X(h(-1),x)$
(cf.  (\ref{Yh(-1)})).

What we really want is a rigorous expression that ``equals'' the
unrigorous expression $X(u,e^{y}x_2)X(v,x_2)$; this will considerably
generalize our interpretation of the unrigorous expression $h(e^{y}
x_2)h(x_2)$ (cf. (\ref{hex1hex2}) and (\ref{hexhex})).  So we replace
$u$ by $(e^{y}x_2)^{L(0)}u$ and $v$ by $x_2^{L(0)}v$ in (\ref{assoc2})
and using basic techniques we get:
\begin{equation}\label{assoc3}
X(u,e^{y}x_2)X(v,x_2) \sim X(Y(e^{yL(0)}u,e^{y}-1)v,x_2),
\end{equation}
and this right-hand side is still rigorous (and the left-hand side
still unrigorous).  But $Y(e^{yL(0)}u,e^{y}-1)$ is exactly Zhu's
operator $Y[u,y]$ in \cite{Z1}, \cite{Z2}, so that
\begin{equation}\label{assoc4}
X(u,e^{y}x_2)X(v,x_2) \sim X(Y[u,y]v,x_2).
\end{equation}
By Zhu's change-of-variables theorem, $x \mapsto Y[u,x]$ defines a new
vertex operator algebra structure on the same vector space $V$ under
suitable conditions; this theorem was a step in Zhu's
vertex-operator-algebraic proof of the modular-invariance properties
of ``characters'' (cf. the comments surrounding (\ref{chiS}) above).
There have been two subsequent treatments of this change-of-variables
theorem, in \cite{L1} and in \cite{H1}, \cite{H2}; in the latter
works, Huang considerably generalized this result (and removed a
hypothesis of Zhu's) using his geometric analysis of the Virasoro
algebra structure in a vertex operator algebra.  The formal relation
(\ref{assoc4}) generalizes to products of several operators.

So we want to bracket the (rigorous) expressions $X(Y[u,y]v,x)$, which
are the correct natural generalization of the expression
(\ref{hexhex}) above (at least with $y_2 = 1$ in (\ref{hexhex}), but
this restriction is a minor issue since $y_2$ can easily be restored).
Keep in mind that the formal relation (\ref{assoc4}) naturally
implements the formal relation (\ref{zeta-s}) in a foundational way
{}from the viewpoint of vertex operator algebras.

But just as in \cite{FLM}, where the Jacobi identity for vertex
operator algebras was the most natural general formula, here we find
that the most natural thing to do is to formulate and prove a new
Jacobi identity, based on the operators $X(Y[u,y]v,x)$, in the general
setting of an arbitrary vertex (operator) algebra, rather than just to
bracket the operators.  It turns out that delta-function expressions
of the type appearing on the right-hand side of (\ref{Lbarbrackets}),
and that in turn ``explained'' the phenomenon of pure monomials in $m$
(as discussed above), arise naturally in complete generality, and when
we ask for a Jacobi identity rather than just a commutator formula in
general, we find that delta-function expressions of this same type
appear on the {\it left-hand side} as well as the right-hand side.
This is another instance of the naturalness of ``Jacobi identities,''
which have symmetries that commutator formulas do not have.  We state
our result for the operators $X(v,x)$ rather than $X(Y[u,y]v,x)$
(i.e., the case where $u$ is the vacuum vector):

\begin{theo}\label{newjacobithm}
In any vertex (operator) algebra $V$, for $u,v \in V$ we have:
\begin{eqnarray}\label{newjacobiiden}
&x_{0}^{-1}\delta
\left(e^{y_{21}}{\displaystyle\frac{x_{1}}{x_{0}}}\right)X(u, x_{1})X(v,
x_{2})-x_{0}^{-1} \delta
\left(-e^{y_{12}}{\displaystyle\frac{x_{2}}{x_{0}}}\right)X(v, x_{2})X(u,
x_{1})&\nonumber \\ &=x_{2}^{-1} \delta
\left(e^{y_{01}}{\displaystyle\frac{x_{1}}{x_{2}}}\right)X(Y[u, -y_{01}]v,
x_{2}),&
\end{eqnarray}
where
\begin{equation}
y_{21} = \log \left(1-{\displaystyle\frac{x_{2}}{x_{1}}}\right),\;\;
y_{12} = \log \left(1-{\displaystyle\frac{x_{1}}{x_{2}}}\right),\;\;
y_{01} = \log \left(1-{\displaystyle\frac{x_{0}}{x_{1}}}\right).\nonumber
\end{equation}
\end{theo}

If we want the {\it commutator} $[X(u, x_{1}),X(v, x_{2})]$, we simply
extract the coefficient of $x_0^{-1}$ (the formal residue in the
variable $x_0$) on both sides, and it turns out that the resulting
right-hand side can be put into an elegant form.  If we replace $u$
and $v$ by expressions of the shape $Y[u,y]v$, we obtain naturally a
formula that generalizes formula (\ref{Lbarbrackets}) (Theorem
\ref{theoremforLbar}) to arbitrary elements of arbitrary vertex
(operator) algebras.  That is, interesting as they are, the phenomena
that we have been discussing concerning central extensions of Lie
algebras of differential operators form extremely special cases of
general vertex-operator-algebraic phenomena.  The detailed
formulations and proofs, and generalizations, are found in \cite{L2}.
Also, my student Antun Milas has generalized some of these results in
a number of directions.

{\small \sc Department of Mathematics, Rutgers University, Piscataway,
NJ 08854}

{\em E-mail address}: lepowsky@math.rutgers.edu


\begin{thebibliography}{FKRW}

\bibitem[BPZ]{BPZ}
A.~A.~Belavin, A.~M.~Polyakov and A.~B.~Zamolodchikov,
 Infinite conformal symmetries in two-dimensional quantum 
field theory,
 {\em Nucl. Phys.} {\bf B241} (1984), 333--380.

\bibitem[Bl]{Bl}
S.~Bloch,
 Zeta values and differential operators on the circle,
 {\em J. Algebra} {\bf 182} (1996), 476--500.

\bibitem[Bo]{Bo}
R.~E.~Borcherds,
 Vertex algebras, Kac-Moody algebras, and the Monster,
 {\em Proc. Natl. Acad. Sci. USA} {\bf 83} (1986), 3068--3071.

\bibitem[DKM]{DKM}
E.~Date, M.~Kashiwara and T.~Miwa,
 Vertex operators and $\tau$ functions---transformation
groups for soliton equations II,
 {\em Proc. Japan Acad. Ser. A Math. Sci.} {\bf 57} (1981), 
387--392.

\bibitem[DL]{DL}
C.~Dong and J.~Lepowsky,
 {\em Generalized Vertex Algebras and Relative
Vertex Operators},
 Progress in Math., Vol. 112, Birkh\"{a}user,
Boston, 1993.

\bibitem[FKRW]{FKRW}
E.~Frenkel, V.~Kac, A.~Radul and W.~Wang,
 $W_{1+\infty}$ and $W(gl_N)$ with central charge $N$,
 {\em Comm. Math. Physics} {\bf 170} (1995), 337--357.

\bibitem[FHL]{FHL}
I.~B.~Frenkel, Y.-Z.~Huang and J.~Lepowsky,
 On axiomatic approaches to vertex operator algebras 
and modules,
 preprint, 1989;
 {\em Memoirs Amer. Math. Soc.} {\bf 104}, 1993.

\bibitem[FLM]{FLM}
I.~B. Frenkel, J.~Lepowsky and A.~Meurman,
 {\em Vertex Operator Algebras and the Monster},
 Pure and Appl. Math., Vol. 134, Academic Press, Boston, 1988.

\bibitem[G]{G}
P.~ Goddard,
 Meromorphic conformal field theory,
 {\em Infinite Dimensional Lie Algebras and Groups, Advanced
Series in Math. Physics,} Vol. 7, ed. V.~Kac, World Scientific,
Singapore, 1989, 556--587.

\bibitem[Hi]{Hi}
H.~Hida,
 {\em Elementary Theory of $L$-functions and Eisenstein Series},
 London Math. Soc. Student Texts, Vol. 26, Cambridge
University Press, Cambridge, 1993.

\bibitem[H1]{H1}
Y.-Z.~Huang,
 Applications of the geometric interpretation of vertex operator
algebras,
 {\em Proc. 20th International Conference on Differential
Geometric Methods in Theoretical Physics, New York, 1991}, ed.
S.~Catto and A.~Rocha, World Scientific, Singapore, 1992, 333--343.

\bibitem[H2]{H2}
Y.-Z.~Huang, 
 {\em Two-dimensional Conformal Geometry and Vertex
Operator Algebras},
 Progress in Math., Vol. 148, Birkh\"{a}user,
Boston, 1997.

\bibitem[KP]{KP}
V.~Kac and D.~Peterson,
 Spin and wedge representations of infinite-dimensional Lie
algebras and groups,
 {\em Proc. Natl. Acad. Sci. USA} {\bf 78} (1981), 3308--3312.

\bibitem[KR]{KR}
V.~Kac and A.~Radul,
 Quasifinite highest weight modules over the Lie algebra of
differential operators on the circle,
 {\em Comm. Math. Physics} {\bf 157} (1993), 429--457.

\bibitem[L1]{L1}
J.~Lepowsky,
 Remarks on vertex operator algebras and moonshine, 
 {\em Proc. 20th International Conference on Differential
Geometric Methods in Theoretical Physics, New York, 1991}, ed.
S.~Catto and A.~Rocha, World Scientific, Singapore, 1992, 362--370.

\bibitem[L2]{L2}
J.~Lepowsky,
 A ``Jacobi identity'' for vertex operator algebras related 
to zeta-function values,
 to appear.

\bibitem[Li1]{Li1}
H.~Li,
 Local systems of vertex operators, vertex superalgebras and 
modules,
 preprint, 1993;
 {\em J. Pure Appl. Alg.} {\bf 109} (1996), 143--195.

\bibitem[Li2]{Li2}
H.~Li,
 Local systems of twisted vertex operators, vertex
superalgebras and twisted modules,
 {\em Contemporary Math.} {\bf 193} (1996), 203--236.

\bibitem[M]{M}
Y.~Matsuo,
 Free fields and quasi-finite representations of $W_{1+\infty}$,
 {\em Physics Lett. B} {\bf 326} (1994), 95--100.

\bibitem[Z1]{Z1}
Y.~Zhu,
 Vertex operators, elliptic functions and modular forms,
 Ph.D. thesis, Yale University, 1990.

\bibitem[Z2]{Z2}
Y.~Zhu,
 Modular invariance of characters of vertex operator algebras,
 {\em J. Amer. Math. Soc.} {\bf 9} (1996), 237--307.

\end{thebibliography}
\end{document}